\documentclass[11pt]{article}
\usepackage{graphicx}
\usepackage{amsfonts}
\usepackage{amsmath, amscd, amsfonts, amssymb, graphicx, color, enumerate,cite}
\usepackage[bookmarksnumbered, colorlinks, plainpages]{hyperref}
\hypersetup{colorlinks=true,linkcolor=red, anchorcolor=green, citecolor=cyan, urlcolor=red, filecolor=magenta, pdftoolbar=true}
\usepackage{theorem}
\newtheorem{Lem}{Lemma}[section]

\newtheorem{theorem}[Lem]{Theorem}
\newtheorem{prop}[Lem]{Proposition}

\newtheorem{remark}[Lem]{Remark}

\newcommand{\qed}{\hbox{\rule{6pt}{6pt}}}

\begin{document}
\title{Note on constants appearing in refined Young inequalities}
\author{Shigeru Furuichi\footnote{E-mail:furuichi@chs.nihon-u.ac.jp}\\
{\small  Department of Information Science, College of Humanities and Sciences, Nihon University,}\\
{\small 3-25-40, Sakurajyousui, Setagaya-ku, Tokyo, 156-8550, Japan}}
\date{}
\maketitle
{\bf Abstract.}
In this short note, we give the refined Young inequality with Specht's ratio by only elementary and direct calculations. The obtained inequality is better than one previously shown by the author in 2012. In addition, we give a new property of Specht's ratio. These imply an alternative proof of the refined Young ineqaulity shown by author in 2012. We also give a remark on the relation to Kantorovich constant. Finally, we give a proposition for a corresponding general function.
\vspace{3mm}


\vspace{3mm}
\vspace{3mm}


\section{Introduction}
Young inequality, $(1-v) a +v b \geq a^{1-v} b^v$ for $a,b>0$ and $0 \leq v \leq 1$ has been refined with constants which are greater than $1$, in many literatures \cite{Dragomir01,Dragomir02,Furuichi2017,FM2011,FM,FMS,LWZ2015}. In this short note, we treat the refinements of Young inequality on ratio type only, although the refinements of Young inequality on difference type has been also studied  by many researchers.

In \cite{Furuichi2012}, we obtained the refined Young inequality with Specht's ratio 
\begin{equation} \label{ineq00}
\frac{(1-v)+vx}{x^v} \geq S(x^r),\quad (x>0)
\end{equation}
where $S(x)$ is Specht's ratio \cite{Spe}  given by $S(x)\equiv \frac{x^{\frac{1}{x-1}}}{e\log x^{\frac{1}{x-1}}}$  and $r \equiv \min\{v,1-v\}$ for $v \in [0,1]$.

In this short note, we show the refined Young inequality with Specht's ratio and a new property of Specht's ratio.
Then these give an alternative proof of the refined Young inequality (\ref{ineq00}).
Finally we consider the general property for which we showed the property on Specht's ratio and Kantorovich constant. Then we state a sufficient condition that such a general property holds for a general function.

\section{Main results}
\begin{theorem}\label{theorem01}
For $v\in [0,1]$ and $x>0$,
\begin{equation} \label{ineq01_theorem01}
\frac{(1-v)+vx}{x^v} \geq S(x)^r,
\end{equation}
where $r \equiv \min\{v,1-v\}$.
\end{theorem}

{\it Proof}:
We set 
$$
g(x)=\log\left\{ (1-v)+v x\right\}-v \log x -r \left(\frac{\log x}{x-1} -1 - \log \frac{\log x}{x-1}\right)
$$
\begin{itemize}
\item[(I)]
We firstly prove $g(x) \geq 0$ for $0 <  x \leq 1$. We calculate the deivative of $g(x)$ by $x$ as
$$
g'(x)=\frac{h(x)}{x(x-1)^2\left\{(1-v)+v x\right\} \log x},
$$
where
$$
h(x) \equiv v(1-v)(x-1)^3\log x+r \left\{(1-v)+v x\right\} \left\{x\left(\log x\right)^2 +(x-1)^2 -(x^2-1)\log x\right\}.
$$
\begin{itemize}
\item[(i)]
For the case of $0 \leq v \leq \frac{1}{2}$, $h(x)=v h_1(x,v)$ where
$$
h_1(x,v)\equiv  (1-v)(x-1)^3\log x+ \left\{(1-v)+v x\right\} \left\{x\left(\log x\right)^2 +(x-1)^2 -(x^2-1)\log x\right\}.
$$
Then we have for $0 < x \leq 1$
$$
\frac{dh_1(x,v)}{dv} = (x-1)\left\{(x-1)^2-2x(x-1) \log x +x \left(\log x\right)^2\right\} \geq 0
$$
due to Lemma \ref{lem02} below.
Thus we have $h_1(x,v) \geq h_1(x,0)$ with 
$$h_1(x,0)=(x-1)^3\log x+x(\log x)^2+(x-1)^2-(x-1)(x+1)\log x.$$
Since $h_1(1,0) = 0$, we assume $0<x<1$.
Then $h_1(x,0) \geq 0$ is equivalent to
$$
(x-1)^2+\frac{x \log x}{x-1} +\frac{x-1}{\log x} -(x+1) \geq 0
$$
dividing the both sides by $(x-1)\log x \geq 0$ for $x>0$. From the arithmetic-geometric mean ineqaulity, it is sufficient to prove the ineqaulity
$$
(x-1)^2+2\sqrt{x}-(x+1) \geq 0
$$
which is true since $(x-1)^2+2\sqrt{x}-(x+1) =\sqrt{x}(\sqrt{x}-1)^2(\sqrt{x} + 2) \geq 0$.
Thus we have $h(x) \geq 0$ which implies $g'(x) \leq 0$. Therefore we have $g(x) \geq g(1) =0$ for $0<x \leq 1$ and $0 \leq v \leq \frac{1}{2}$.
Note that $\lim_{x\to 1}\log\frac{\log x}{x-1}=0$ by $\lim_{x\to 1} \frac{\log x}{x-1}=1$.
\item[(ii)]
For the case of $\frac{1}{2}\leq v \leq 1$ and $0<x \leq 1$, $h(x)=(1-v) h_2(x,v)$ where
$$
h_2(x,v)\equiv v(x-1)^3\log x+\left\{(1-v)+v x\right\} \left\{x(\log x)^2+(x-1)^2 -(x^2-1)\log x\right\}.
$$
Then we have
$$
\frac{dh_2(x,v)}{dv} =(x-1)k_2(x)
$$
where
$$
k_2(x)\equiv (x-1)^2-2(x-1)\log x +x(\log x)^2.
$$
Since $k_2'(x)=2\left(\sqrt{x}-\frac{1}{\sqrt{x}}\right)^2 +(\log x)^2 \geq 0$, we have $k_2(x) \leq k_2(1) = 0$ for $0< x \leq 1$. Thus we have $\frac{dh_2(x,v)}{dv} \geq 0$ for $0<x \leq 1$ which implies $h_2(x,v) \geq h_2(x,\frac{1}{2})$.
So we here prove $2 h_2(x,\frac{1}{2}) =(x-1)^3 \log x+(x+1) \left\{ x(\log x)^2 +(x-1)^2-(x^2-1)\log x\right\}\geq 0$ which is equivalent to
$$
(x-1)^2+\frac{x(x+1)\log x}{x-1} +\frac{(x+1)(x-1)}{\log x} \geq (x+1)^2
$$
dividing the both sides by $(x-1)\log x \geq 0$ for $x > 0$.
From the arithmetic-geometric mean inequality, it is sufficient to prove the inequality
$$
(x-1)^2+2\sqrt{x}(x+1) -(x+1)^2 \geq 0
$$
which is simplified as $x+1 \geq 2 \sqrt{x}$ and this is trivial.
Thus we have $h(x) \geq 0$ which implies $g'(x) \leq 0$. Therefore we have $g(x) \geq g(1) =0$ for $0<x \leq 1$ and $\frac{1}{2} \leq v \leq 1$.
\end{itemize}
\item[(II)]
We finally prove $g(x) \geq 0$ for the case of $x \geq 1$.
From Lemma \ref{lem01} below, we see
$g(x) \geq f(x)$
where
$$
f(x) =\log\left\{ (1-v)+v x\right\}-v \log x -r \left(\frac{2}{x+1} -1 - \log \frac{2}{x+1}\right).
$$
Then we have
$$
f'(x) = \frac{(x-1) f_1(x)}{x(x+1)^2\left\{(1-v) + v x\right\}}
$$
where
$$
f_1(x)=v(1-v-r)x^2+(1-v)(2v-r)x+v(1-v).
$$
Since $f_1'(x) = 2v(1-v-r)x+(1-v)(2v -r)$,
we have $f_1'(x) = 2v(1-2v)x +v(1-v) \geq 0$ for $0 \leq v \leq \frac{1}{2}$ and $x \geq 1$, also have
$f_1'(x)=(1-v)(3v-1) \geq 0$ for $\frac{1}{2} \leq v \leq 1$ and $x \geq 1$.
Thus we have $f_1(x) \geq f_1(1) = 4v(1-v)-r \geq 0$. Thus we have $f'(x) \geq 0$ and we then have $f(x) \geq f(1) =0$ for $x \geq 1$.
\end{itemize}

\hfill \qed

\begin{Lem}\label{lem02}
For $0<x \leq 1$, we have
\begin{equation}\label{ineq01_lem02}
(x-1)^2+x (\log x)^2 \leq 2x(x-1)\log x
\end{equation}
\end{Lem}
{\it Proof:}
Since $(x-1)\log x \geq 0$ for $0 < x \leq 1$, the ineqaulity (\ref{ineq01_lem02}) is equivalent to the inequality
$$
\frac{x-1}{\log x}+\frac{x\log x}{x-1} \geq 2 x.
$$
By the arithmetic-geometric mean inequality, we have the first inequality in the following
$$
\frac{x-1}{\log x}+\frac{x\log x}{x-1} \geq 2 \sqrt{x} \geq 2 x.
$$

\hfill \qed

\begin{Lem}\label{lem01}
For $x \geq 1$, we have
$$\frac{\log x}{x-1} -\log \frac{\log x}{x-1} \leq \frac{2}{x+1} -\log  \frac{2}{x+1}.$$
\end{Lem}
{\it Proof:}
Since the function $y-\log y$ is decreasing for $0< y \leq 1$,
$\frac{2}{x+1} \leq \frac{\log x}{x-1} \leq 1$ implies the desired result.

\hfill \qed

\begin{remark}
Since Kantorovich constant $K(x) \equiv\frac{(x+1)^2}{4x}$ is geater than or equal to $S(x)$, we have $K(x)^r \geq S(x)^r$ for $0 \leq r \leq 1$ so that we easily obtain the inequality
\begin{equation}\label{ineq00_Remark2.4}
\frac{(1-v)+xv}{x^v} \geq K(x)^r \geq S(x)^r,\quad 0\leq r \leq \frac{1}{2}
\end{equation}
from the result in \cite{ZSF2011}, where $r \equiv \min\{v,1-v\}$. For the convenience to the readers, we give the proof of the first inequality in (\ref{ineq00_Remark2.4}). Since $r=v$ for the case $v \in [0,1/2]$, we consider the function $w_1(v,x)\equiv (1-v)+v x-\left(\frac{x+1}{2}\right)^{2v}$ for $v \in [0,1/2]$ and $x>0$. By simple calculations, we got $\frac{d^2w_1(v,x)}{dv^2} = -4^{1-v} (x+1)^{2v}\left( \log\left(\frac{x+1}{2}\right)\right)^2 \leq 0$ with $w_1(0,x)=w_1(1/2,x)=0$, thus we have $w_1(v,x) \geq 0$   for $v \in [0,1/2]$ and $x>0$. Similarly, since $r=1-v$ for the case $v \in [1/2,1]$, we consider the function $w_2(v,x)\equiv \frac{(1-v)+v x}{x} -\left(\frac{x+1}{2x}\right)^{2(1-v)}$ for $v \in [1/2,1]$ and $x>0$. Since we have 
$\frac{d^2w_2(v,x)}{dv^2} = -4^v\left(\frac{x+1}{x}\right)^{2(1-v)}\left(\log\left(\frac{x+1}{2x}\right)\right)^2 \leq 0$ with $w_2(1/2,x)=w_2(1,x)=0$, we have $w_2(v,x) \geq 0$ for $v \in [1/2,1]$ and $x>0$. 

We also have $K(x^r) \geq S(x^r)$ for $x>0$ and $0\leq r \leq 1$. For the convenience of the readers, we give a proof of the inequality $K(x) \geq S(x)$ for $x>0$, which is equivalent to
$w_1(x) \geq 0$, where $w_1(x) \equiv 2\log(x+1)-2\log 2-\log x-\frac{\log x}{x-1} +1+\log\left(\frac{\log x}{x-1}\right)$.
Since $\frac{\log t}{t-1} \geq \frac{2}{t+1}$ for $t>0$, we have only to prove the inequality $w_2(x) \geq 0$, where $w_2(x) \equiv 2\log(x+1)-2\log 2-\log x-\frac{\log x}{x-1} +1+\log\left(\frac{2}{x+1}\right)=1-\frac{\log x}{x-1}+\log\left(\frac{x+1}{2x}\right)$. Then we have
$w_2'(x)=\frac{2(1-x)+(x+1)\log x}{(x+1)(x-1)^2}$. We easily find 
$w_2'(x) < 0$ for $0<x < 1$ and $w_2'(x) > 0$ for $x > 1$, since we use again the relation $\frac{\log t}{t-1} \geq \frac{2}{t+1}$ for $t>0$.
Thus we have $w_2(x) > w_2(1) = 1$ for $x >0$ with $x \neq 1$. It is trivial $w_2(1)=K(1)=S(1)=1$ for $x=1$.
\end{remark}

\begin{remark}
For $x>0$ and $0 \leq r \leq 1$, we have the inequality
$$
K(x)^r \geq K(x^r).
$$
The proof is done by the following. The above inequality is equivalent to the inequality
$
g_r(x) \geq 0,
$
where
$$
g_r(x)\equiv r\log(x+1)-\log(x^r+1)-r\log2+\log 2.
$$
Since $g_r'(x)=\frac{r(1-x^{1-r})}{(x+1)(x^r+1)}$, we have $g_r'(1)=0$, $g_r'(x) \leq 0$ for $0<x \leq 1$ and $g_r'(x) \geq 0$ for $x \geq 1$.
Thus we have $g_r(x) \geq g_r(1)=0$.
\end{remark}

\begin{remark}\label{Remark2.6}
We have no ordering between
$S(x)^r$ and  $K(x^r)$, since we have the following examples.
When $x=5$ and $r=0.2$, $S(x)^r-K(x^r) \simeq 0.0384244$. When $x=30$ and $r=0.49$,  $S(x)^r-K(x^r) \simeq -0.0162285$.
When $x=0.1$ and $r=0.4$, $S(x)^r-K(x^r) \simeq 0.0534669$.
When $x=0.01$ and $r=0.49$, $S(x)^r-K(x^r) \simeq -0.0955003$.
\end{remark}
\begin{theorem} \label{theorem02}
We have the  inequality  
$$
S(x)^r \geq S(x^r)
$$
for  $x>0$ and $0 \leq r \leq 1$.
\end{theorem}

To prove Theorem \ref{theorem02}, we prepare the following lemma.
\begin{Lem}\label{lemma1.8}
Put $k(t) \equiv -(t-1)^3+3t(t-1)(\log t)^2 -t(t+1) (\log t)^3.$
Then $k(t) \geq 0$ for $0<t \leq 1$ and $k(t) \leq 0$ for $t \geq 1$.
\end{Lem}
{\it Proof}:
It is trivial that $k(1)=0$. We firstly consider the case of $t \geq 1$.
By calculations, we have $k^{(3)}(t)=-3t^{-2} l(t)$, where $l(t)\equiv 2(t-1)^2+2(t-1)\log t+(2t-1)(\log t)^2$. Then we have $l'(t)=-2-2t^{-1}+4t+2(3-t^{-1})\log t+2(\log t)^2$ and $l''(t) =2t^{-2}\left\{ t(2t+3) +(2t+1)\log t\right\} \geq 0$ for $t \geq 1$. Thus we have $l'(t) \geq l'(1) =0$ and then $l(t) \geq l(1)=0$, namely $k^{(3)}(t) \leq 0$ for $t \geq 1$. Thus we have $k''(t) \leq k''(1)=0$, where $k''(t) = t^{-1}\left\{-6(t-1)^2+12(t-1)\log t -3(t+1)(\log t)^2-2t(\log t)^3\right\}$. And then $k'(t) \leq k'(1) =0$, where
 $k'(t)=-3(t-1)^2+6(t-1)\log t+3(t-2)(\log t)^2-(2t+1)(\log t)^3$. Therefore we have $k(t) \leq k(1)=0$.

Since we have just shown that
$$
-(t-1)^3+3t(t-1)(\log t)^2-t(t+1)(\log t)^3 \leq 0,\quad (t \geq1),
$$
we have  
$$
-\left(\frac{1}{s}-1\right)^3+\frac{3}{s}\left(\frac{1}{s}-1\right)(\log s)^2+\frac{1}{s}\left(\frac{1}{s}+1\right)(\log s)^3 \leq 0, \quad (0<s \leq 1)
$$
putting $t=s^{-1}$. Thus we have $k(s) \geq 0$ for $0<s \leq 1$, by multipying $s^3$ to both sides in the above inequality.

\hfill
\qed

{\it Proof of Theorem \ref{theorem02}}:
We set 
\begin{eqnarray*}
u(x,r) &\equiv& r \log(S(x)) -\log(S(x^r)) \\
&=& r\left(\frac{\log x}{x-1}-1-\log\frac{\log x}{x-1}\right)
-\left(\frac{r\log x}{x^r-1}-1-\log\frac{r\log x}{x^r-1}\right).
\end{eqnarray*}
Then we have for any $x >0$ and $0\leq r \leq 1$
$$\frac{d^2u(x,r)}{dr^2} =\frac{-(x^r-1)^3+3x^r(x^r-1)(\log x^r)^2-x^r(x^r+1)(\log x^r)^3}{r^2(x^r-1)^3} \leq 0$$ 
by Lemma \ref{lemma1.8} with $t=x^r$. With $u(x,1) =0$ and $\lim_{r\to 0} u(x,r) =0$ by $\lim_{r\to 0} \frac{x^r-1}{r} = \log x$, we have $u(x,r) \geq 0$ for any  $x >0$ and $0\leq r \leq 1$. Since it is known that $S(x) \geq 1$, $u(x,r) \geq 0$ implies the inequality $S(x)^r \geq S(x^r)$ for any $x >0$ and $0\leq r \leq 1$.

\hfill
\qed

Thus Theorem \ref{theorem01} and Theorem \ref{theorem02} gives an alternative proof of the inequality given in (\ref{ineq00}). In addition, the inequality in (\ref{ineq01_theorem01}) gives a better bound than one in (\ref{ineq00}).

\section{Concluding remarks}
Summarizing the ordering the constants appeared in this short note, we have the following relations for $r=\min\left\{v,1-v\right\}$, $v \in [0,1]$ and $x>0$.
\[\begin{array}{l}
\frac{{\left( {1 - v} \right) + vx}}{{{x^v}}} \ge K{\left( x \right)^r} \ge K\left( {{x^r}} \right)\\
\,\,\,\,\,\,\,\,\,\,\,\,\,\,\,\,\,\,\,\,\,\,\,\,\,\,\,\,\,\,\,\,\,\, \vee  \!\!\;\!\!\;|\,\,\,\,\,\,\,\,\,\,\,\,\,\,\,\, \vee \!\!\;\!\!\;\!\!\;\!\!\;\!\!\;\!\!\;|\,\,\,\\
\frac{{\left( {1 - v} \right) + vx}}{{{x^v}}} \ge S{\left( x \right)^r} \ge S\left( {{x^r}} \right)
\end{array}\]
In addition, we have no ordering between $S(x)^r$ and $K(x^r)$ as stated in Remark \ref{Remark2.6}.

It is known that Kantorovich constant and  Specht's ratio have similar properties:
\begin{itemize}
\item They take infinity in the limits as $x \to 0$ and $x \to \infty$.
\item They take $1$ when $x=1$.
\item They are monotone decreasing functions on $(0,1)$.
\item They are monotone increasing functions on $(1,\infty)$.
\end{itemize}
As we have seen in Section 2, properties $K(x)^r \geq K(x^r)$ and $S(x)^r \geq S(x^r)$   for $x>0$ and $0 \leq r \leq 1$, make us to consider the following general result. To state our proposition, we define a geometrically convex function $f:J\to (0,\infty)$ where $J \subset (0,\infty)$. If the function $f:J\to (0,\infty)$ where $J \subset (0,\infty)$ satisfies the inequality $f(a^{1-v}b^v) \leq f(a)^{1-v}f(b)^v$ for $a,b >0$ and $v\in [0,1]$, then the function $f$ is called a geometrically convex (multiplicatively convex) function.  See \cite{Nicul2000} e.g., for the  geometrically convex (multiplicatively convex) function.
We easily find that if the geometrically convex function $f:J\to (0,\infty)$ is decreasing, then $f$ is also (usual) convex function.

It is known that the function $f(x)$ is a geometrically convex if and only if the function $\log f(e^x)$ is (usual) convex. See e.g., \cite{Audenaert_LAA2013}.

\begin{prop}\label{prop3.1}
Assume the function $f: (0,\infty) \to (0,\infty)$ satisfies $ f(1) \leq 1$. If $f$ is a geometrically convex function,
 then we have $f(x)^r \geq f(x^r)$ for $x>0$ and $0 \leq r \leq 1$.
\end{prop}

\noindent{\it Proof}.
Put $a=1$, $b=x$ and $v=r$ in the definition of a geometrically convex function, $f(a^{1-v}b^v) \leq f(a)^{1-v}f(b)^v$. Then we obtain the inequality $f(x)^r \geq f(x^r)$ for $x>0$ and $0 \leq r \leq 1$, since we have the condition $0< f(1) \leq 1$.

\hfill \qed

\begin{remark}
If we assume $f$ is a twice differentiable function, we have the following alternative proof for Proposition \ref{prop3.1}.
We put $g(r,x) = r \log(f(x)) -\log(f(x^r))$. 
It is known in \cite[Proposition 4.3]{Nicul2000} that the geometrically convexity (multiplicatively convexity) of $f$  is equivalent to the following condition \eqref{con_01},  under the assumption that $f$ is twice differentiable. 
\begin{equation} \label{con_01}
D(x) := f(x)f'(x) +x \left(f(x)f''(x)-f'(x)^2\right) \geq 0,\,\,\,\, (x>0)
\end{equation}

Then  we calculate
$$
\frac{d^2g(r,x)}{dr^2}=-\frac{x^r(\log x)^2}{f(x^r)^2} 
\left\{f(x^r)\left(f'(x^r)+x^rf''(x^r)\right)-x^rf'(x^r)^2\right\} \leq 0,
$$
by the condition (\ref{con_01}).
It is easy to check $g(0,x)\geq 0,\,\, g(1,x)=0$ so that we have $g(r,x) \geq 0$ for any $x>0$ and $0\leq r \leq 1$.
We thus obtain the inequality $f(x)^r \geq f(x^r)$ for $x>0$ and $0 \leq r \leq 1$.
\end{remark}

Thus we can make a judgement whether the function $f$ is a geometrically convex or not, by using $D(x)$ or $\frac{d^2\log f(e^x)}{dx^2}$.
We can find $D(x) = \frac{(x+1)^2}{8x^2} \geq 0$ for Kantorovich constant $K(x)$. The expression of $D(x) $ is complicated for Specht's ratio $S(x)$. But we could check that it takes non-negative values for $x>0$  by the numerical computations.  In addition, for the function $agr_v(x) \equiv \frac{(1-v) + vx}{x^v}$ defined for $x>0$ with $0\leq v \leq 1$,
we find that $D(x) = v(1-v)x^{-2v} \geq 0$ so that $arg(x)^r \geq arg(x^r)$ for $0 \leq r \leq 1$ and $x>0$.

\begin{remark}
As for the converse in Proposition \ref{prop3.1}, under the assumption we do not impose the condition $0<f(1)\leq 1$,
the claim that the inequality $f(x)^r \geq f(x^r)$ for $x>0,\,\,0\leq r\leq 1$,  implies $f$ is a geometrically convex function, is not true in general. 
The function $f(x):=\frac{1-x}{1+x}$ on $(0,1)$ is not a geometrically convex since $D(x)=-\frac{2(x^2+1)}{(x+1)^4} < 0$, (also $\frac{d^2\log f(e^x)}{dx^2}=\frac{-2e^x(e^{2x}+1)}{(e^{2x}-1)^2} < 0$). Actually, the function $f(x)$ is a geometrically concave. On the other hand, we have the inequality $f(x)^r \geq f(x^r)$ for $0< x < 1$ and $0\leq r \le 1$. Indeed, we have
\begin{equation}\label{remark3.3_ineq01}
\left(\frac{1-x}{1+x}\right)^r \geq \frac{1-x}{1+x} \geq \frac{1-x^r}{1+x^r}.
\end{equation}
The first inequality is due to $0< \frac{1-x}{1+x} < 1$ and the second one is easy calculation as $(1-x)(1+x^r)-(1+x)(1-x^r)=2(x^r-x) \geq 0$.

However, we have not found any examples such that the function $f$ satisfies $0<f(1)\leq 1$ and the inequality $f(x)^r \geq f(x^r)$, but it is not a geometrically convex. In the above example $f(x)=\frac{1-x}{1+x}$ on $(0,1)$, we see actually $f(1):=\mathop {\lim }\limits_{x \nearrow 1} f\left( x \right) =0$. We considered the open interval $(0,1)$ not $(0,1]$ in the above example, to avoid the case $0^0$ in the left hand side of the inequalities \eqref{remark3.3_ineq01}. We also have not found the proof that the inequality $f(x)^r \geq f(x^r)$ for the function $f:(0,\infty)\to(0,\infty)$ under the assumption $0< f(1) \leq 1$ , implies the geometrically convexity of $f$.
\end{remark}

Note that Kantorovich constant is convex by $K''(x) = \frac{x^3}{2} \geq 0$ but Specht's ratio is not convex in general, since $S''(10) \simeq -0.00156565$ by numerical computations. If we impose the log-convexity on Proposition \ref{prop3.1}, then the inequality (\ref{con_01}) holds only for $x \geq 1$. However the Kantorovich constant is not log-convex. These show that it may be difficult to simplify the condition (\ref{con_01}), more than this. Note here that  if the function $f:J\to (0,\infty)$ satisfies  the inequality $f((1-v)a+vb)\leq f(a)^{1-v}f(b)^b$ for $a,b>0$ and $v\in [0,1]$, then $f$ is called log-convex function.  It is easy to find that all log-convex functions are also (usual) convex function by the use of the weighted arithmetic-geometric mean inequality.

In the papers \cite{Furuichi2017,FM}, we proved the refined Young inequalities for a special case such as $0<x \leq 1$:
$$
\frac{(1-v)+vx}{x^v} \geq \frac{1}{1-\frac{v(1-v)(x-1)^2}{2}},
\quad \frac{(1-v)+vx}{x^v} \geq 1+\frac{2^vv(1-v)(x-1)^2}{(x+1)^{v+1}}
$$
However, the inequality $f(x)^r \geq f(x^r)$ for $0< x \leq 1$ does not hold in general  for these lower bounds of $\frac{(1-v)+vx}{x^v}$.
Actually, at least the sufficient condition given in Proposition \ref{prop3.1} such that
$D(x) \geq 0$ for $0<x \leq 1$ does not hold in general for these lower bounds of $\frac{(1-v)+vx}{x^v}$. They can be checked by numerical computations.

We conclude this short note by stating that a property $f(x)^r \geq f(x^r)$ has been studied in \cite{Wada2014,Wada2018} in details by the different approach. Such a property is named there as power monotonicity and applied to some known results in the framework of operator theory.

\section*{Acknowledgement}
The author was partially supported by JSPS KAKENHI Grant Number 16K05257.

\end{document}